# A Further Study of an $L^2$-norm Based Test for the Equality of Several Covariance Functions

Jia Guo and Jin-Ting Zhang

National University of Singapore

**Abstract**

For the multi-sample equal covariance function (ECF) testing problem, Zhang (2013) proposed an $L^2$-norm based test. However, its asymptotic power and finite sample performance have not been studied. In this paper, its asymptotic power is investigated under some mild conditions. It is shown that the $L^2$–norm based test is root-$n$ consistent. In addition, intensive simulation studies demonstrate that in terms of size-controlling and power, the $L^2$-norm based test outperforms the dimension-reduction based test proposed by Fremdt et al. (2013) when the functional data are less correlated or when the effective signal information is located in high frequencies. Two real data applications are also presented to demonstrate the good performance of the $L^2$-norm based test.

Short Title: Testing Equality of Several Covariance Functions

## 1 Introduction

In the recent decades, functional data, which are easily recorded in the form of curves or images with the rapid development of data collecting technologies, attract much attention in the statistical literature. The early efforts were devoted to the description of functional data. A good survey is given by Ramsay and Silverman (2005). In the recent decade, much work has Jia Guo (E-mail: jia.guo@u.nus.edu) is PhD candidate, Jin-Ting Zhang (E-mail: stazjt@nus.edu.sg) is Associate Professor, Department of Statistics and Applied Probability, National University of Singapore, Singapore 117546.
The work was financially supported by the National University of Singapore Academic Research grant R-155-000-164-112.

been done in hypothesis testing about the mean functions of one or several functional populations. The reader is referred to Zhang (2013) and references therein. For testing the equality of the mean functions of two functional populations, Ramsay and Silverman (2005) described a pointwise t-test. This pointwise t-test is a natural extension of the classical t-test but it has some obvious drawbacks: it is conducted at any individual time point, and it does not give an overall conclusion. To overcome this problem, Zhang et al. (2010b) proposed an $L^2$-norm based test. For the functional one-way ANOVA problem, Ramsay and Silverman (2005) suggested a pointwise $F$-test which is an extension of the classical $F$-test to the context of functional data analysis. This pointwise $F$-test has the same limitation as the pointwise t-test. Some alternatives are then proposed. Cuevas et al. (2004) proposed an $L^2$-norm based functional ANOVA test. However, their test adopts time-consuming Monte Carlo simulations to approximate the null distribution. Zhang and Liang (2013) studied a so-called GPF test, whose test statistic was obtained via globalizing the pointwise $F$-test for the functional one-way ANOVA problem. In the context of functional linear models, Zhang and Chen (2007) studied an $L^2$-norm based test for a general linear hypothesis testing problem and Shen and Faraway (2004) considered a functional $F$-test to compare two nested functional linear models. Via intensive simulations, Górecki and Smaga (2015) presented an exhaustive comparison of a number of existing functional hypothesis testing procedures and concluded that the GPF test works quite well for the functional one-way ANOVA problem.

Besides testing the equality of the mean functions of one or several functional populations, some novel and effective methods have also been proposed for testing the equality of two covariance functions. For example, for testing the equality of the covariance functions of two functional populations, Zhang and Sun (2010) proposed an $L^2$-norm based test while Fremdt et al. (2013) studied a dimension-reduction based test which is an extension of the work of Panaretos et al. (2010) to the non-Gaussian case. The testing procedure of Fremdt et al. (2013) was obtained via projecting the observations onto a suitably chosen finite-dimensional space. However, little work has been done for testing the equality of several covariance functions. We may refer to this problem as the multi-sample equal covariance function (ECF) testing problem. This multi-sample ECF testing problem is encountered frequently in many areas. For example, in the functional one-way and two-way ANOVA, we usually assume that the covariance functions of different samples are the same. However, in real data analysis, this assumption may not be true and a formal test may be needed before applying the previously-mentioned testing procedures for the functional one-way or two-way ANOVA. For this multi-sample ECF testing problem, Zhang (2013) (Ch. 10) described an $L^2$-norm based test, which is simple to implement and easy to interpret. However, its asymptotic power has not been studied. In addition, no simulation results are given to demonstrate its finite-sample performance. In this paper, we present a further study on this $L^2$-norm based test via studying its asymptotic power. As a result, we show that it is a root-n consistent test. We also demonstrate its finite-sample performance via comparing it with Fremdt et al. (2013)'s dimension-reduction based test through intensive simulations. We found that when the functional data are less correlated or when the effective signal information is located in high frequencies, the $L^2$-norm based test is more powerful than the afore-mentioned dimension-reduction based test.

The rest of the paper is organized as follows. Section 2 presents a simple review of the $L^2$-norm based test. Section 3 studies its asymptotic power. Two simulation studies and two real data applications are presented in Sections 4 and 5 respectively. The technical proofs of the main results are given in the Appendix.

## 2 The $L^2$-norm Based Test

In this section, we give a brief review of the $L^2$-norm based test proposed in Zhang (2013) for the $k$-sample equal-covariance function (ECF) testing problem.

The ECF testing problem is defined as below. Let $y_{i1}(t), y_{i2}(t), \cdots, y_{in_i}(t),\ i = 1, 2, \cdots, k$ be $k$ independent functional samples over a given finite time interval $\mathcal{T} = [a, b]$, $-\infty < a < b < \infty$, which satisfy

$$\begin{aligned} y_{ij}(t) &= \eta_i(t) + v_{ij}(t),\ j = 1, 2, \cdots, n_i, \\ v_{i1}(t), v_{i2}(t), \cdots, v_{in_i}(t) &\stackrel{i.i.d.}{\sim} \mathrm{SP}(0, \gamma_i);\ i = 1, 2, \cdots, k, \end{aligned} \tag{2.1}$$

where $\eta_1(t), \eta_2(t), \cdots, \eta_k(t)$ model the unknown group mean functions of the $k$ samples, $v_{ij}(t),\ j = 1, 2, \cdots, n_i,\ i = 1, 2, \cdots, k$ denote the subject-effect functions which follow a stochastic process with mean function 0 and covariance function $\gamma_i(s, t),\ i = 1, 2, \cdots, k$ respectively. Throughout this paper, we assume that $\mathrm{tr}(\gamma_i) < \infty$ and $\eta_i(t) \in \mathcal{L}^2(\mathcal{T}),\ i = 1, 2, \cdots, k$, where $\mathcal{L}^2(\mathcal{T})$ denotes the Hilbert space formed by all the squared integrable functions over $\mathcal{T}$ with the inner-product defined as $< f, g >= \int_{\mathcal{T}} f(t)g(t)dt,\ f,\ g \in \mathcal{L}^2(\mathcal{T})$. It is often of interest to test the equality of the $k$ covariance functions:

$$H_0 : \gamma_1(s,t) \equiv \gamma_2(s,t) \equiv \cdots \equiv \gamma_k(s,t), \quad \text{for all } s,t \in \mathcal{T}. \tag{2.2}$$

Based on the given $k$ functional samples (2.1), the group mean functions $\eta_i(t),\ i = 1,2,\cdots,k$ and the covariance functions $\gamma_i(s,t),\ i = 1,2,\cdots,k$ can be unbiasedly estimated as

$$\begin{aligned} \hat{\eta}_i(t) &= \bar{y}_i(t) = n_i^{-1} \textstyle\sum_{j=1}^{n_i} y_{ij}(t),\ i = 1,2,\cdots,k, \\ \hat{\gamma}_i(s,t) &= (n_i-1)^{-1} \textstyle\sum_{j=1}^{n_i} [y_{ij}(s) - \bar{y}_i(s)][y_{ij}(t) - \bar{y}_i(t)],\ i = 1,2,\cdots,k. \end{aligned} \tag{2.3}$$

It is easy to show that $\hat{\gamma}_i(s,t),\ i = 1,2,\cdots,k$ are independent and $\mathrm{E}\hat{\gamma}_i(s,t) = \gamma_i(s,t),\ i = 1,2,\cdots,k$. Further, the estimated subject-effect functions can be written as

$$\hat{v}_{ij}(t) = y_{ij}(t) - \bar{y}_i(t),\ j = 1,2,\cdots,n_i;\ i = 1,2,\cdots,k. \tag{2.4}$$

When the null hypothesis (2.2) holds, we let $\gamma(s,t)$ denote the common covariance function of the $k$ samples. It can be estimated by the following pooled sample covariance function

$$\hat{\gamma}(s,t) = \sum_{i=1}^{k} (n_i - 1)\hat{\gamma}_i(s,t)/(n-k), \tag{2.5}$$

where $\hat{\gamma}_i(s,t),\ i = 1,2,\cdots,k$ are given in (2.3) and throughout, $n = \sum_{i=1}^{k} n_i$ denotes the total sample size.

For further investigation, the following assumptions are imposed.

Assumption A

1. The $k$ samples (2.1) are Gaussian.
2. As $n \to \infty$, the $k$ sample sizes satisfy $n_i/n \to \tau_i \in (0,1),\ i = 1,2,\cdots,k$.
3. The variance functions $\gamma_i(s,t), i = 1,2,\cdots,k$ are uniformly bounded.

It is easy to note that $\hat{\gamma}_i(s,t) - \hat{\gamma}(s,t)$ measures the difference between the $i$-th sample covariance function (2.3) and the pooled sample covariance function (2.5), which should be small when the null hypothesis holds. Based on this, Zhang (2013) proposed the following so-called $L^2$-norm based test statistic for the $k$-sample ECF testing problem (2.2):

$$T_n = \sum_{i=1}^{k} (n_i - 1) \int_{\mathcal{T}} \int_{\mathcal{T}} [\hat{\gamma}_i(s,t) - \hat{\gamma}(s,t)]^2 ds dt, \tag{2.6}$$

which summarizes all the squared differences between the $k$ sample covariance functions and the pooled sample covariance function. Therefore, when the null hypothesis holds, $T_n$ will be small and otherwise large.

Lemma **??** in the Appendix states that the test statistic $T_n$ is asymptotically a $\chi^2$-type mixture. Therefore, the null distribution of $T_n$ can be approximated by the well-known Welch-Satterthwaite $\chi^2$-approximation. By this method, Zhang (2013) approximated the null distribution of $T_n$ using that of a random variable

$$R \sim \beta\chi_d^2. \tag{2.7}$$

The parameters $\beta$ and $d$ are determined via matching the first two moments of $T_n$ and $R$, which are given by

$$\beta = \frac{\mathrm{tr}(\varpi^{\otimes 2})}{\mathrm{tr}(\varpi)},\ d = (k-1)\kappa, \kappa = \frac{\mathrm{tr}^2(\varpi)}{\mathrm{tr}(\varpi^{\otimes 2})}, \tag{2.8}$$

where $\varpi[(s_1, t_1), (s_2, t_2)]$ denotes the covariance function of $\sqrt{n-k}[\hat{\gamma}(s,t) - \gamma(s,t)]$, $\mathrm{tr}(\varpi) = \int_{\mathcal{T}}\int_{\mathcal{T}} \varpi\,[(s,t),(s,t)]\,dsdt$ and $\mathrm{tr}(\varpi^{\otimes 2}) = \int_{\mathcal{T}}\int_{\mathcal{T}}\int_{\mathcal{T}}\int_{\mathcal{T}} \varpi^2\,[(s_1,t_1),(s_2,t_2)]\,ds_1dt_1ds_2dt_2$. Under the Gaussian assumption A1, it is easy to verify that

$$\begin{aligned} \varpi[(s_1,t_1),(s_2,t_2)] &= \gamma(s_1,s_2)\gamma(t_1,t_2) + \gamma(s_1,t_2)\gamma(s_2,t_1), \\ \mathrm{tr}(\varpi) = \mathrm{tr}^2(\gamma) + \mathrm{tr}(\gamma^{\otimes 2}), &\ \ \mathrm{tr}(\varpi^{\otimes 2}) = 2\mathrm{tr}^2(\gamma^{\otimes 2}) + 2\mathrm{tr}(\gamma^{\otimes 4}), \end{aligned} \tag{2.9}$$

where $\mathrm{tr}(\gamma) = \int_{\mathcal{T}} \gamma(t,t)dt$, $\mathrm{tr}(\gamma^{\otimes 2}) = \int_{\mathcal{T}}\int_{\mathcal{T}} \gamma^2(s,t)dsdt$ and

$$\mathrm{tr}(\gamma^{\otimes 4}) = \int_{\mathcal{T}}\int_{\mathcal{T}}\int_{\mathcal{T}}\int_{\mathcal{T}} \gamma(t,u_1)\gamma(u_1,u_2)\gamma(u_2,u_3)\gamma(u_3,t)du_1du_2du_3dt.$$

To conduct the $L^2$-norm based test, we need to estimate the parameters $\beta$ and $d$ based on the data. There are two methods for estimating the parameters $\beta$ and $\kappa$, one is the naive method, and the other is the bias-reduced method. Let $\hat{\beta}$ and $\hat{\kappa}$ denote the estimators of $\beta$ and $\kappa$. The naive estimators of $\beta$ and $\kappa$ are obtained via replacing $\mathrm{tr}(\varpi), \mathrm{tr}^2(\varpi)$ and $\mathrm{tr}(\varpi^{\otimes 2})$ in (2.8) respectively with their naive estimators $\mathrm{tr}(\hat{\varpi}), \mathrm{tr}^2(\hat{\varpi})$ and $\mathrm{tr}(\hat{\varpi}^{\otimes 2})$:

$$\hat{\beta} = \frac{\mathrm{tr}(\hat{\varpi}^{\otimes 2})}{\mathrm{tr}(\hat{\varpi})},\ \ \hat{\kappa} = \frac{\mathrm{tr}^2(\hat{\varpi})}{\mathrm{tr}(\hat{\varpi}^{\otimes 2})}, \tag{2.10}$$

where under the Gaussian assumption A1 and based on (2.9), we have

$$\begin{aligned} \hat{\varpi}[(s_1,t_1),(s_2,t_2)] &= \hat{\gamma}(s_1,s_2)\hat{\gamma}(t_1,t_2) + \hat{\gamma}(s_1,t_2)\hat{\gamma}(s_2,t_1), \\ \mathrm{tr}(\hat{\varpi}) = \mathrm{tr}^2(\hat{\gamma}) + \mathrm{tr}(\hat{\gamma}^{\otimes 2}), &\ \ \mathrm{tr}(\hat{\varpi}^{\otimes 2}) = 2\mathrm{tr}^2(\hat{\gamma}^{\otimes 2}) + 2\mathrm{tr}(\hat{\gamma}^{\otimes 4}). \end{aligned} \tag{2.11}$$

The bias-reduced estimators of $\beta$ and $\kappa$ are obtained via replacing $\mathrm{tr}(\varpi), \mathrm{tr}^2(\varpi)$ and $tr(\varpi^{\otimes 2})$ in (2.8) respectively with their bias-reduced estimators $\widehat{\mathrm{tr}(\varpi)}$, $\widehat{\mathrm{tr}^2(\varpi)}$ and $\widehat{\mathrm{tr}(\varpi^{\otimes 2})}$:

$$\hat{\beta} = \frac{\widehat{\mathrm{tr}(\varpi^{\otimes 2})}}{\widehat{\mathrm{tr}(\varpi)}}, \quad \hat{\kappa} = \frac{\widehat{\mathrm{tr}^2(\varpi)}}{\widehat{\mathrm{tr}(\varpi^{\otimes 2})}} \tag{2.12}$$

where under the Gaussian assumption A1 and based on (2.9), we have

$$\begin{aligned} \widehat{\mathrm{tr}(\varpi)} = \widehat{\mathrm{tr}^2(\gamma)} + \widehat{\mathrm{tr}(\gamma^{\otimes 2})}, \quad \widehat{\mathrm{tr}^2(\varpi)} = \left[\widehat{\mathrm{tr}(\varpi)}\right]^2, \\ \widehat{\mathrm{tr}(\varpi^{\otimes 2})} = 2\left[\widehat{\mathrm{tr}(\gamma^{\otimes 2})}\right]^2 + 2\mathrm{tr}(\hat{\gamma}^{\otimes 4}), \end{aligned} \tag{2.13}$$

with

$$\begin{aligned} \widehat{\mathrm{tr}^2(\gamma)} &= \tfrac{(n-k)(n-k+1)}{(n-k-1)(n-k+2)}\left[\mathrm{tr}^2(\hat{\gamma}) - \tfrac{2\mathrm{tr}(\hat{\gamma}^{\otimes 2})}{n-k+1}\right], \\ \widehat{\mathrm{tr}(\gamma^{\otimes 2})} &= \tfrac{(n-k)^2}{(n-k-1)(n-k+2)}\left[\mathrm{tr}(\hat{\gamma}^{\otimes 2}) - \tfrac{\mathrm{tr}^2(\hat{\gamma})}{n-k}\right]. \end{aligned} \tag{2.14}$$

Note that under the Gaussian assumption A1, $\widehat{\mathrm{tr}^2(\gamma)}$ and $\widehat{\mathrm{tr}(\gamma^{\otimes 2})}$ are the unbiased estimators of $\mathrm{tr}^2(\gamma)$ and $\mathrm{tr}(\gamma^{\otimes 2})$ respectively and when the data are not Gaussian, they may be asymptotically unbiased under some further assumptions. Notice also that in the expression (2.13), the unbiased estimator of $\mathrm{tr}(\gamma^{\otimes 4})$ is not incorporated since it is quite challenging to obtain a simple and useful unbiased estimator of $\mathrm{tr}(\gamma^{\otimes 4})$.

The following theorem shows that under some mild conditions, the estimators, $\hat{\beta}$ and $\hat{\kappa}$, of $\beta$ and $\kappa$ are consistent.

**Theorem 1.** *Under Assumptions A1∼A3 and the null hypothesis (2.2), as $n \to \infty$, we have $\hat{\beta} \xrightarrow{p} \beta$, $\hat{\kappa} \xrightarrow{p} \kappa$ for both the naive and bias-reduced methods. In addition, $\hat{T}_n(\alpha) \xrightarrow{p} \tilde{T}_0(\alpha)$, where $\hat{T}_n(\alpha) = \hat{\beta}\chi^2_{(k-1)\hat{\kappa}}(\alpha)$ is the estimated critical value of $T_n$ and $\tilde{T}_0(\alpha) = \beta\chi^2_{(k-1)\kappa}(\alpha)$ is its approximate theoretical critical value.*

By some simple algebra, we have $\beta < \lambda_{\max} < \infty$ and $\kappa \leq m$ where $\lambda_{\max}$ is the largest eigenvalue of $\varpi\left[(s_1, t_1), (s_2, t_2)\right]$ and $m$ is the number of all the positive eigenvalues. Then it is easy to verify that $\tilde{T}_0(\alpha) < \infty$ when $m$ is a finite number.

However, under the null hypothesis, when the sample sizes $n_i, i = 1, 2, \cdots, k$ of the $k$ samples (2.1) are small, Theorem 1 is no longer valid so that the Welch-Satterthwaite $\chi^2$-approximation is also no longer applicable. To overcome this difficulty, a random permutation method is proposed to approximate the critical values of $T_n$. This method can also be used when the data are non-Gaussian. The random permutation method can be described as follows.

Firstly, we randomly reorder the pooled estimated subject-effect functions (2.4) so that a random permutation sample $\hat{v}_l^*(t),\ l = 1, 2, \cdots, n$ is obtained where $n$ is the total sample size as defined before. We then use the first $n_1$ permuted subject-effect functions to form the first permutation sample $\hat{v}_{1j}^*(t),\ j = 1, 2, \cdots, n_1$, the next $n_2$ permuted subject-effect functions to form the second permutation sample $\hat{v}_{2j}^*(t),\ j = 1, 2, \cdots, n_2$ and so on. The permutation test statistic $T_n^*$ is computed similar to the computation of the original $L^2$-norm based test statistic $T_n$ as described in (2.6) but now based on the $k$ permuted functional samples. That is,

$$T_n^* = \sum_{i=1}^k (n_i - 1) \int_{\mathcal{T}} \int_{\mathcal{T}} \left[\hat{\gamma}_i^*(s,t) - \hat{\gamma}^*(s,t)\right]^2 dsdt,$$

where

$$\hat{\gamma}_i^*(s,t) = (n_i - 1)^{-1} \textstyle\sum_{j=1}^{n_i} \hat{v}_{ij}^*(s)\hat{v}_{ij}^*(t), i = 1, 2, \cdots, k,$$
$$\hat{\gamma}^*(s,t) = \textstyle\sum_{i=1}^k (n_i - 1)\hat{\gamma}_i^*(s,t)/(n-k).$$

Repeating the above process a large number of times, we can get a sample of $T_n^*$ and use the sample upper $100\alpha$-percentile $\hat{T}_n^*(\alpha)$ to estimate the critical value of $T_n^*$. Using this critical value, we then conduct the associated random permutation test. If $T_n > \hat{T}_n^*(\alpha)$, we reject the null hypothesis (2.2).

The following theorem shows that under the null hypothesis the permutation test statistic $T_n^*$ converges in distribution to the same limit test statistic $T_0$ of $T_n$ where $T_0$ is defined in Lemma **??** and hence $\hat{T}_n^*(\alpha)$ will also tend to $T_0(\alpha)$ in distribution as $n \to \infty$. Thus the size of the permutation test tends to the nominal size.

**Theorem 2.** *Under Assumptions A1∼A3 and the null hypothesis (2.2), as $n \to \infty$, we have $T_n^* \xrightarrow{d} T_0$ and $\hat{T}_n^*(\alpha) \xrightarrow{d} T_0(\alpha)$. Hence, the size of the random permutation test $P(T_n > \hat{T}_n^*(\alpha)) \to P(T_n > T_0(\alpha))$ where $\hat{T}_n^*(\alpha)$ is the estimated upper $100\alpha$-percentile of $T_n^*$ based on the permutation samples and $T_0(\alpha)$ is the theoretical critical value of $T_n$.*

## 3 Asymptotic Power of the $L^2$-norm Test

Zhang (2013) did not study the asymptotic power of the $L^2$-norm based test $T_n$. In this section, we study its asymptotic power under the following local alternative:

$$H_1 : \gamma_i(s,t) = \gamma(s,t) + (n_i - 1)^{-1/2} d_i(s,t),\ i = 1, 2, \cdots, k, \tag{3.1}$$

where $d_1(s,t), d_2(s,t), \cdots, d_k(s,t)$ are some fixed bivariate functions, independent of $n$ and $\gamma(s,t)$ is some covariance function.

For further study, we can re-write the $L^2$-norm based test statistic $T_n$ (2.6) as

$$T_n = \int_{\mathcal{T}} \int_{\mathcal{T}} \text{SSB}_n(s,t) ds dt, \tag{3.2}$$

where

$$\text{SSB}_n(s,t) = \sum_{i=1}^{k} (n_i - 1)[\hat{\gamma}_i(s,t) - \hat{\gamma}(s,t)]^2, \tag{3.3}$$

which summarizes the squared differences between the individual sample covariance functions $\hat{\gamma}_i(s,t), i = 1, 2, \cdots, k$ and the pooled sample covariance function $\hat{\gamma}(s,t)$ for any given $(s,t) \in \mathcal{T}^2$.

Before we state the main results, we give an alternative expression of $\text{SSB}(s,t)$ which is helpful for deriving the asymptotic power of $T_n$. For any $s, t \in \mathcal{T}$, $\text{SSB}_n(s,t)$ can be expressed as

$$SSB_n(s,t) = \mathbf{z}_n(s,t)^T[\mathbf{I}_k - \mathbf{b}_n\mathbf{b}_n^T/(n-k)]\mathbf{z}_n(s,t) = \boldsymbol{z}_n(s,t)^T \boldsymbol{W}_n \boldsymbol{z}_n(s,t), \tag{3.4}$$

where

$$\mathbf{z}_n(s,t) = [z_1(s,t), z_2(s,t), \cdots, z_k(s,t)]^T, \quad \boldsymbol{W}_n = \mathbf{I}_k - \mathbf{b}_n\mathbf{b}_n^T/(n-k), \tag{3.5}$$

with

$$\begin{aligned} z_i[s,t] &= \sqrt{n_i - 1}[\hat{\gamma}_i(s,t) - \gamma(s,t)], \ i = 1, 2, \cdots, k, \\ \mathbf{b}_n &= [\sqrt{n_1 - 1}, \sqrt{n_2 - 1}, \cdots, \sqrt{n_k - 1}]^T. \end{aligned}$$

Since $\mathbf{b}_n^T\mathbf{b}_n/(n-k) = 1$, it is easy to verify that $\boldsymbol{W}_n$ is an idempotent matrix with rank $k-1$. In addition, as $n \to \infty$, we have

$$\boldsymbol{W}_n \to \boldsymbol{W} := \mathbf{I}_k - \mathbf{b}\mathbf{b}^T, \text{with } \mathbf{b} = [\sqrt{\tau_1}, \sqrt{\tau_2}, \cdots, \sqrt{\tau_k}]^T, \tag{3.6}$$

where $\tau_i, i = 1, 2, \cdots, k$ are given in Assumption A2. Note that $\boldsymbol{W}$ in (3.6) is also an idempotent matrix of rank $k-1$, which has the following singular value decomposition:

$$\boldsymbol{W} = \mathbf{U} \begin{pmatrix} \mathbf{I}_{k-1} & \mathbf{0} \\ \mathbf{0}^T & 0 \end{pmatrix} \mathbf{U}^T, \tag{3.7}$$

where the columns of $\mathbf{U}$ are the eigenvectors of $\boldsymbol{W}$.

Let $\tilde{\boldsymbol{d}}(s,t) = [\mathbf{I}_{k-1}, \mathbf{0}]\mathbf{U}^T\mathbf{d}(s,t)$ where $\mathbf{d}(s,t) = [d_1(s,t), d_2(s,t), \cdots, d_k(s,t)]^T$ with $d_i(s,t), i = 1, 2, \cdots, k$ given in (3.1). Let $\lambda_r, r = 1, 2, \cdots$ be the eigenvalues of $\varpi[(s_1,t_1),(s_2,t_2)]$ with only the first $m$ eigenvalues being positive and $\phi_r(s,t)$, $r = 1, 2, \cdots$, be the associated eigenfunctions. Define

$$\delta_r^2 = ||\int_{\mathcal{T}}\int_{\mathcal{T}} \tilde{\boldsymbol{d}}(s,t)\phi_r(s,t)dsdt||^2, r = 1, 2, \cdots, \tag{3.8}$$

which measure the information of $\tilde{\boldsymbol{d}}(s,t)$ projected on the eigenfunctions $\phi_r(s,t), r = 1, 2, \cdots$, of $\varpi[(s_1,t_1),(s_2,t_2)]$. Theorem 3 below gives the asymptotic distribution of $T_n$ under the local alternative (3.1).

**Theorem 3.** *Under Assumptions A1∼A3 and the local alternative (3.1), as* $n \to \infty$, *we have* $T_n \xrightarrow{d} T_1$ *with*

$$T_1 \stackrel{d}{=} \Sigma_{r=1}^m \lambda_r A_r + \Sigma_{r=m+1}^{\infty} \delta_r^2,$$

*where* $A_r \sim \chi^2_{k-1}(\lambda_r^{-1}\delta_r^2)$, $r = 1, 2, \cdots, m$, *are independent and* $\delta_r^2, r = m+1, m+2, \cdots, \infty$ *are defined in (3.8).*

Theorem 4 below shows that under the local alternative (3.1), $T_n$ is asymptotically normal. Theorem 5 below shows that the $L^2$-norm based test can detect the local alternative (3.1) with probability 1 provided that the information provided by $\mathbf{d}(s,t)$ diverges. That is, the $L^2$-norm based test is root-$n$ consistent. In both Theorems 4 and 5, the quantities $\delta_r^2, r = 1, 2, \cdots$ have been defined in (3.8). Let these quantities satisfy the following condition:

$$\max_r \delta_r^2 \to \infty. \tag{3.9}$$

This condition describes a situation when the information projected onto at least one eigenfunction tends to $\infty$.

**Theorem 4.** *Under Assumptions A1∼A3, the local alternative (3.1), and condition (3.9), as* $n \to \infty$, *we have* $\frac{T_n - E(T_n)}{\sqrt{Var(T_n)}} \xrightarrow{d} N(0,1)$.

**Theorem 5.** *Under Assumptions A1∼A3, the local alternative (3.1), and condition (3.9), as* $n \to \infty$, *the proposed* $L^2$*-norm based test has asymptotic power 1, i.e.,* $P(T_n > \hat{T}_n(\alpha)) \to 1$, *where* $\hat{T}_n(\alpha)$ *is the estimated critical value of* $T_n$ *defined in Theorem 1.*

We now study the consistency property of the random permutation test. Theorem 6 shows that the random permutation test is also root-$n$ consistent.

**Theorem 6.** *Under Assumptions A1∼A3, the local alternative (3.1), and condition (3.9), as* $n \to \infty$, *the power* $P(T_n > \hat{T}_n^*(\alpha)) \to 1$.

## 4 Simulation Studies

In Section 2, we described three methods for approximating the null distribution of the $L^2$-norm based test: a naive method, a bias-reduced method, and a random-permutation method. The associated $L^2$-norm based tests may be denoted as $L^2_{nv}, L^2_{br}$ and $L^2_{rp}$ respectively. Recently, Fremdt et al. (2013) described two dimension-reduction methods for testing the equality of the covariance functions of two functional samples. Their first test can be applied to both Gaussian and non-Gaussian functional data while the second one can only be used for Gaussian functional data. For convenience, we refer to these two tests as $FHK_D$ and $FHK_G$ respectively. In this section, we shall present two simulations. In Simulation 1, we shall compare the performances of $L^2_{nv}, L^2_{br}$ and $L^2_{rp}$ and in Simulation 2, we shall compare $L^2_{nv}, L^2_{br}$ and $L^2_{rp}$ against $FHK_D$ and $FHK_G$.

### 4.1 Data generating

In the simulations, for $i = 1, 2, \cdots, k$, the $i$-th functional sample will be generated from the following model:

$$\begin{aligned} &y_{ij}(t) = \eta_i(t) + v_{ij}(t),\ \eta_i(t) = \mathbf{c}_i^T[1, t, t^2, t^3]^T,\ v_{ij}(t) = \mathbf{b}_{ij}^T \mathbf{\Psi}_i(t),\ t \in [0, 1], \\ &\mathbf{b}_{ij} = [b_{ij1}, b_{ij2}, \cdots, b_{ijq}]^T,\ b_{ijr} \stackrel{d}{=} \sqrt{\lambda_r} z_{ijr},\ r = 1, 2, \cdots, q; j = 1, 2, \cdots, n_i, \end{aligned} \tag{4.1}$$

where the parameter vectors $\mathbf{c}_i = [c_{i1}, c_{i2}, c_{i3}, c_{i4}]^T$ for the group mean function $\eta_i(t)$ can be flexibly specified, the random variables $z_{ijr},\ r = 1, 2, \cdots, q$ are i.i.d. with mean 0 and variance 1, $\mathbf{\Psi}_i(t) = [\psi_{i1}(t), \psi_{i2}(t), \cdots, \psi_{iq}(t)]^T$ is a vector of $q$ basis functions and the variance components $\lambda_r,\ r = 1, 2, \cdots, q$ are positive and decreasing in $r$, and the number of the basis functions, $q$, is an odd positive integer. These tuning parameters help specify the group mean functions $\eta_i(t) = c_{i1} + c_{i2}t + c_{i3}t^2 + c_{i4}t^3$ and the covariance function $\gamma_i(s, t) = \mathbf{\Psi}_i(s)^T diag(\lambda_1, \lambda_2, \cdots, \lambda_q)\mathbf{\Psi}_i(t) = \sum_{r=1}^q \lambda_r \psi_{ir}(s)\psi_{ir}(t),\ i = 1, 2, \cdots, k$. For simplicity, we assume that the design time points for all the functions $y_{ij}(t),\ j = 1, 2, \cdots, n_i,\ i = 1, 2, \cdots, k$ are the same and are specified as $t_j = (j-1)/(J-1),\ j = 1, 2, \cdots, J$, where $J$ is some positive integer. In practice, these functions can be observed at different design time points. In this case, some smoothing technique, such as those discussed in Zhang and Chen (2007), Zhang et al. (2010a) can be used to reconstruct the functions $y_{ij}(t),\ j = 1, 2, \cdots, n_i,\ i = 1, 2, \cdots, k$ and then to evaluate them at a common grid of time points. The latter simulation setup will be time-consuming to carry out and we did not explore it in this simulation.

We now specify the parameters in (4.1). To specify the group mean functions $\eta_1(t), \eta_2(t), \cdots, \eta_k(t)$, we set $\mathbf{c}_1 = [1, 2.3, 3.4, 1.5]^T$ and $\mathbf{c}_i = \mathbf{c}_1 + (i-1)\delta\mathbf{u},\ i = 1, 2, \cdots, k$, where the constant vector $\mathbf{u}$ specifies the direction of these differences. We set $\delta = 0.1$ and $\mathbf{u} = [1, 2, 3, 4]^T/\sqrt{30}$ which is a unit vector. Then we specify the covariance function $\gamma_i(s,t)$. For simplicity, we set $\lambda_r = a\rho^{r-1},\ r = 1, 2, \cdots, q$, for some $a > 0$ and $0 < \rho < 1$. Notice that the tuning parameter $\rho$ not only determines the decay rate of $\lambda_1, \lambda_2, \cdots, \lambda_q$, but also determines how the simulated functional data are correlated: when $\rho$ is close to 0, $\lambda_1, \lambda_2, \cdots, \lambda_q$ will decay very fast, indicating that the simulated functional data are highly correlated; and when $\rho$ is close to 1, $\lambda_r,\ r = 1, 2, \cdots, q$ will decay slowly, indicating that the simulated functional data are nearly uncorrelated. To define the basis functions $\mathbf{\Psi}_i(t)$, we firstly generate a vector of $q$ basis functions $\boldsymbol{\phi}(t) = [\phi_1(t), \phi_2(t), \cdots, \phi_q(t)]^T$ and we select $\phi_1(t) = 1,\ \phi_{2r}(t) = \sqrt{2}sin(2\pi rt),\ \phi_{2r+1}(t) = \sqrt{2}cos(2\pi rt),\ t \in [0, 1],\ r = 1, 2, \cdots, (q-1)/2$. Then we specify our basis functions $\mathbf{\Psi}_i(t)$ via the following relationship: $\psi_{ir}(t) = \phi_r(t),\ r = 1, 3, 4, \cdots, q$ but $\psi_{i2}(t) = \phi_2(t) + (i-1)\omega, i = 1, 2, \cdots, k$. That is, we obtain $k$ different bases via shifting the second basis function of the $i$-th basis with $(i-1)\omega$ steps. This allows that the differences of the $k$ covariance functions $\gamma_i(s,t), i = 1, 2, \cdots, k$ are controlled by the tuning parameter $\omega$ since we actually have

$$\gamma_i(s,t) = \gamma_1(s,t) + (i-1)\lambda_2(\phi_2(s) + \phi_2(t))\omega + (i-1)^2\lambda_2\omega^2,\ i = 1, 2, \cdots, k. \quad (4.2)$$

Further, we set $a = 1.5,\ q = 11$ and the number of design time points $J = 180$. Finally, we specify two cases of the distribution of the i.i.d. random variables $z_{ijr},\ r = 1, 2, \cdots, q;\ j = 1, 2, \cdots, n_i;\ i = 1, 2, \cdots, k$: $z_{ijr} \overset{i.i.d.}{\sim} N(0, 1)$ and $z_{ijr} \overset{i.i.d.}{\sim} t_4/\sqrt{2}$, allowing to generate Gaussian and non-Gaussian functional data respectively with $z_{ijr}$ having mean 0 and variance 1. Notice that the $t_4/\sqrt{2}$ distribution is chosen since it has nearly heaviest tails among the $t$-distributions with finite first two moments.

## 4.2 Simulation 1: a comparison of $L^2_{nv}, L^2_{br}$ and $L^2_{rp}$

In this simulation, to check the finite sample performance of $L^2_{nv}, L^2_{br}$ and $L^2_{rp}$, we let the number of groups $k = 5$. We consider three cases of the sample size vector: $\mathbf{n}_1 = [20, 25, 22, 18, 16]$, $\mathbf{n}_2 = [35, 30, 40, 32, 38]$ and $\mathbf{n}_3 = [80, 75, 85, 82, 70]$, representing the small, moderate, and large sample size cases. We also consider four correlation cases, i.e., $\rho = 0.1, 0.3, 0.5, 0.7$, representing the highly, moderately, less correlated, and nearly independent situations. For given model configurations, the required functional samples are generated. The p-values of $L^2_{nv}, L^2_{br}$ and $L^2_{rp}$ are then computed. Notice that the p-values of $L^2_{rp}$ is obtained via 500 runs of permutations. We reject the null hypothesis if the calculated p-values are smaller than the nominal significance level $\alpha = 5\%$. We repeat the above simulation process 10000 times to get the empirical sizes or powers of $L^2_{nv}, L^2_{br}$ and $L^2_{rp}$.

Table 1 shows the empirical sizes and powers (in percentages) of $L^2_{nv}, L^2_{br}$ and $L^2_{rp}$ when $z_{ijr},\ r = 1, 2, \cdots, q;\ j = 1, 2, \cdots, n_i;\ i = 1, 2, \cdots, k \overset{i.i.d.}{\sim} N(0, 1)$. We have the following conclusions:

- In terms of size controlling, $L^2_{nv}$ works well when the functional data are highly correlated but it becomes rather conservative (with the empirical size to be as small as 4.02%) when the correlation of functional data is reduced. $L^2_{br}$ generally works well for various settings and it becomes better with increasing the sample sizes. $L^2_{rp}$ is quite liberal (with the empirical size to be as large as 8.40%), especially when the functional data are nearly independent. However, it performs better with increasing the sample sizes.
- In terms of powers, $L^2_{br}$ is comparable or have higher powers than $L^2_{nv}$ and $L^2_{rp}$ when their empirical sizes are comparable.
- Overall, when the functional data are Gaussian, $L^2_{br}$ outperforms $L^2_{nv}$ and $L^2_{rp}$.

Note that the inflated sizes of $L^2_{rp}$ may be due to the small number of runs of permutations, which is 500. However, increasing this number requires much more computational efforts and we did not adopt this strategy for time saving.

Table 2 shows the empirical sizes and powers of $L^2_{nv}, L^2_{br}$ and $L^2_{rp}$ when $z_{ijr},\ r = 1, 2, \cdots, q;\ j = 1, 2, \cdots, n_i;\ i = 1, 2, \cdots, k \overset{i.i.d.}{\sim} t_4/\sqrt{2}$, representing the cases when the functional data are non-Gaussian. We have the following conclusions:

- In terms of size controlling, both $L^2_{nv}$ and $L^2_{br}$ do not work since their empirical sizes are too large compared with the nominal size 5%. This is expected since the formulas (2.11) used for computing the approximated null distributions are based on the Gaussian assumption A1.

Table 1: Empirical sizes and powers (in percentages) of $L^2_{nv}, L^2_{br}$ and $L^2_{rp}$ for Simulation 1 when $z_{ijr},\ r=1,\cdots,q;\ j=1,\cdots,n_i;\ i=1,\cdots,k \overset{i.i.d.}{\sim} N(0,1)$.

| | $\mathbf{n}_1$=[20,25,22,18,16] | | | | $\mathbf{n}_2$=[35,30,40,32,38] | | | | $\mathbf{n}_3$=[80,75,85,82,70] | | | |
|---|---|---|---|---|---|---|---|---|---|---|---|---|
| $\rho$ | $\omega$ | $L^2_{nv}$ | $L^2_{br}$ | $L^2_{rp}$ | $\omega$ | $L^2_{nv}$ | $L^2_{br}$ | $L^2_{rp}$ | $\omega$ | $L^2_{nv}$ | $L^2_{br}$ | $L^2_{rp}$ |
| | 0.00 | 4.29 | 4.55 | 5.93 | 0.00 | 4.90 | 5.09 | 5.75 | 0.00 | 5.14 | 5.23 | 5.51 |
| | 0.70 | 22.50 | 23.14 | 24.75 | 0.60 | 24.71 | 25.14 | 25.27 | 0.40 | 17.78 | 17.97 | 18.41 |
| 0.1 | 1.00 | 50.50 | 51.17 | 50.81 | 0.80 | 52.39 | 52.90 | 51.76 | 0.60 | 54.61 | 54.87 | 53.81 |
| | 1.50 | 87.62 | 88.09 | 85.22 | 1.10 | 87.31 | 87.66 | 85.35 | 0.75 | 83.86 | 83.98 | 83.09 |
| | 2.20 | 99.13 | 99.17 | 98.31 | 1.50 | 99.36 | 99.39 | 98.68 | 0.90 | 96.75 | 96.80 | 96.53 |
| | 0.00 | 4.28 | 4.69 | 5.81 | 0.00 | 5.12 | 5.37 | 6.17 | 0.00 | 5.03 | 5.13 | 5.26 |
| | 0.36 | 22.66 | 23.45 | 25.08 | 0.30 | 24.39 | 24.95 | 25.81 | 0.20 | 22.84 | 23.09 | 23.08 |
| 0.3 | 0.55 | 51.97 | 53.05 | 52.11 | 0.45 | 59.96 | 60.57 | 58.98 | 0.30 | 56.46 | 56.79 | 56.78 |
| | 0.80 | 83.84 | 84.49 | 82.06 | 0.56 | 82.30 | 82.86 | 80.87 | 0.36 | 78.02 | 78.27 | 77.31 |
| | 1.20 | 98.96 | 99.07 | 97.68 | 0.85 | 99.47 | 99.49 | 99.04 | 0.45 | 95.24 | 95.34 | 94.67 |
| | 0.00 | 4.48 | 5.07 | 6.80 | 0.00 | 4.71 | 5.03 | 5.79 | 0.00 | 4.82 | 5.01 | 5.09 |
| | 0.30 | 29.19 | 30.57 | 32.01 | 0.20 | 23.18 | 24.06 | 25.02 | 0.15 | 29.81 | 30.22 | 30.30 |
| 0.5 | 0.40 | 50.05 | 51.72 | 52.44 | 0.30 | 52.02 | 52.96 | 52.35 | 0.20 | 54.57 | 55.15 | 54.38 |
| | 0.60 | 84.06 | 84.85 | 82.66 | 0.40 | 80.31 | 80.99 | 78.95 | 0.25 | 77.07 | 77.45 | 76.50 |
| | 1.00 | 99.48 | 99.52 | 98.65 | 0.80 | 99.98 | 99.98 | 99.94 | 0.40 | 99.66 | 99.67 | 99.61 |
| | 0.00 | 4.02 | 5.11 | 8.40 | 0.00 | 4.39 | 5.11 | 6.47 | 0.00 | 5.30 | 5.69 | 6.27 |
| | 0.22 | 23.31 | 25.62 | 30.15 | 0.16 | 21.45 | 22.99 | 24.78 | 0.11 | 23.88 | 24.65 | 25.37 |
| 0.7 | 0.35 | 53.48 | 55.70 | 57.52 | 0.25 | 53.19 | 54.88 | 54.97 | 0.17 | 58.64 | 59.43 | 58.97 |
| | 0.54 | 87.95 | 88.87 | 86.92 | 0.35 | 83.87 | 84.81 | 83.05 | 0.21 | 79.87 | 80.46 | 79.36 |
| | 0.80 | 99.06 | 99.18 | 97.97 | 0.50 | 98.74 | 98.82 | 97.94 | 0.30 | 98.54 | 98.58 | 98.19 |

- The performance of $L^2_{rp}$, on the other hand, is comparable with those cases presented in Table 1. That is, its empirical sizes are liberal when the functional data are less correlated but they are getting better with increasing the sample sizes.
- Thus, when the functional data are not Gaussian, $L^2_{rp}$ may work for large samples but $L^2_{nv}$ and $L^2_{br}$ do not work at all.

Table 2: Empirical sizes and powers (in percentages) of $L^2_{nv}, L^2_{br}$ and $L^2_{rp}$ for Simulation 1 when $z_{ijr},\ r=1,\cdots,q;\ j=1,\cdots,n_i;\ i=1,\cdots,k \overset{i.i.d.}{\sim} t_4/\sqrt{2}$.

| | $\mathbf{n}_1$=[20,25,22,18,16] | | | | $\mathbf{n}_2$=[35,30,40,32,38] | | | | $\mathbf{n}_3$=[80,75,85,82,70] | | | |
|---|---|---|---|---|---|---|---|---|---|---|---|---|
| $\rho$ | $\omega$ | $L^2_{nv}$ | $L^2_{br}$ | $L^2_{rp}$ | $\omega$ | $L^2_{nv}$ | $L^2_{br}$ | $L^2_{rp}$ | $\omega$ | $L^2_{nv}$ | $L^2_{br}$ | $L^2_{rp}$ |
| 0.1 | 0.00 | 41.60 | 42.38 | 7.00 | 0.00 | 49.38 | 49.75 | 6.04 | 0.00 | 58.60 | 58.75 | 5.56 |
| | 0.90 | 60.65 | 61.36 | 29.03 | 0.72 | 67.06 | 67.46 | 23.14 | 0.60 | 79.88 | 79.93 | 26.71 |
| | 1.50 | 86.84 | 87.35 | 60.64 | 1.20 | 92.37 | 92.49 | 60.85 | 0.85 | 94.02 | 94.08 | 58.20 |
| | 2.50 | 98.71 | 98.76 | 86.42 | 1.80 | 99.38 | 99.40 | 86.62 | 1.20 | 99.58 | 99.60 | 87.38 |
| | 3.50 | 99.74 | 99.75 | 92.08 | 2.80 | 99.99 | 99.99 | 94.97 | 1.50 | 99.97 | 99.97 | 95.32 |
| 0.3 | 0.00 | 48.17 | 49.45 | 7.19 | 0.00 | 55.26 | 55.87 | 6.18 | 0.00 | 63.98 | 64.29 | 5.13 |
| | 0.45 | 62.34 | 63.24 | 25.14 | 0.38 | 74.47 | 74.90 | 24.68 | 0.28 | 84.48 | 84.60 | 22.31 |
| | 0.80 | 86.68 | 87.29 | 58.06 | 0.60 | 91.28 | 91.44 | 54.96 | 0.46 | 96.88 | 96.92 | 60.58 |
| | 1.70 | 99.45 | 99.47 | 90.09 | 1.00 | 99.40 | 99.41 | 85.52 | 0.65 | 99.77 | 99.79 | 86.74 |
| | 3.20 | 99.98 | 99.98 | 94.77 | 2.00 | 100.00 | 100.00 | 96.20 | 0.75 | 99.93 | 99.93 | 92.72 |
| 0.5 | 0.00 | 56.78 | 58.41 | 8.10 | 0.00 | 65.23 | 66.03 | 6.86 | 0.00 | 75.09 | 75.47 | 5.73 |
| | 0.40 | 74.60 | 75.74 | 34.75 | 0.30 | 83.47 | 84.01 | 28.39 | 0.20 | 91.46 | 91.69 | 24.40 |
| | 0.70 | 92.57 | 92.99 | 66.00 | 0.40 | 91.68 | 91.95 | 46.23 | 0.30 | 97.61 | 97.69 | 52.39 |
| | 1.00 | 98.31 | 98.43 | 82.65 | 0.80 | 99.69 | 99.70 | 87.56 | 0.50 | 99.90 | 99.91 | 88.68 |
| | 2.00 | 99.97 | 99.97 | 93.95 | 1.50 | 99.99 | 99.99 | 95.81 | 0.60 | 100.00 | 100.00 | 94.62 |
| 0.7 | 0.00 | 68.67 | 71.37 | 9.94 | 0.00 | 77.79 | 79.51 | 7.22 | 0.00 | 86.81 | 87.34 | 6.12 |
| | 0.24 | 76.59 | 78.43 | 23.97 | 0.25 | 90.91 | 91.62 | 28.31 | 0.15 | 95.25 | 95.49 | 19.74 |
| | 0.52 | 92.17 | 92.88 | 59.40 | 0.40 | 97.15 | 97.33 | 58.57 | 0.25 | 99.25 | 99.32 | 52.32 |
| | 0.88 | 98.91 | 99.01 | 83.59 | 0.70 | 99.89 | 99.90 | 87.88 | 0.40 | 99.96 | 99.96 | 87.05 |
| | 1.80 | 99.99 | 100.00 | 94.45 | 1.20 | 100.00 | 100.00 | 95.97 | 0.55 | 100.00 | 100.00 | 96.55 |

### 4.3 Simulation 2: a comparison of $L^2_{nv}, L^2_{br}$ and $L^2_{rp}$ against $FHK_D$ and $FHK_G$

In this simulation, we shall use the simulation codes kindly provided by Dr. Fremdt via email communication. To compare $L^2_{nv}, L^2_{br}$ and $L^2_{rp}$ against $FHK_D$ and $FHK_G$, we set the number of groups $k=2$ and consider the sample size $\mathbf{n}_1=[25,22]$, $\mathbf{n}_2=[30,40]$ and $\mathbf{n}_3=[75,85]$. We also specify $\rho=0.1,0.3,0.5,0.7$ to consider the four cases when the functional samples have high, moderate, low, very low correlations. Since smoothing is needed to conduct $FHK_D$ and $FHK_G$, we choose 49 Fourier basis functions to smooth the simulated functions. Notice that $FHK_D$ and $FHK_G$ require selecting $d$, the number of empirical functional principal components. To avoid this time-consuming selection, we instead just consider: $d = 1, 2, 3, 4$, hoping that the important signals in functional data are located at low principal components. Actually, we shall use the same method as described in the beginning of this section to generate the functional samples for the simulations presented in this subsection so that the main differences between the two covariance functions are located at the first two basis functions as indicated in (4.2).

Table 3: Empirical sizes (in percentages) of $FHK_D, FHK_G, L^2_{nv}, L^2_{br}$ and $L^2_{rp}$ for Simulation 2 when $z_{ijr},\ r = 1, \cdots, q;\ j = 1, \cdots, n_i;\ i = 1, 2 \overset{i.i.d.}{\sim} N(0, 1)$.

| | | $FHK_D$ | | | | $FHK_G$ | | | | $L^2_{nv}$ | $L^2_{br}$ | $L^2_{rp}$ |
|---|---|---|---|---|---|---|---|---|---|---|---|---|
| **n** | $\rho$ | $d=1$ | $d=2$ | $d=3$ | $d=4$ | $d=1$ | $d=2$ | $d=3$ | $d=4$ | | | |
| $\mathbf{n}_1$ | 0.1 | 4.87 | 3.95 | 3.58 | 3.64 | 5.17 | 4.98 | 4.82 | 4.19 | 4.24 | 4.63 | 5.47 |
| | 0.3 | 5.50 | 4.40 | 3.57 | 3.43 | 5.55 | 5.37 | 4.55 | 4.58 | 3.81 | 4.20 | 5.56 |
| | 0.5 | 4.71 | 4.38 | 3.44 | 3.13 | 5.30 | 5.10 | 4.74 | 4.73 | 3.64 | 4.29 | 5.53 |
| | 0.7 | 4.79 | 4.25 | 4.02 | 3.46 | 5.30 | 5.32 | 5.14 | 4.89 | 2.90 | 4.20 | 6.35 |
| $\mathbf{n}_2$ | 0.1 | 5.34 | 6.31 | 6.47 | 9.56 | 5.08 | 4.91 | 4.51 | 4.62 | 4.42 | 4.65 | 5.21 |
| | 0.3 | 5.73 | 6.07 | 6.39 | 9.03 | 5.54 | 5.43 | 4.86 | 4.97 | 4.43 | 4.78 | 5.35 |
| | 0.5 | 5.50 | 5.62 | 6.79 | 8.92 | 5.29 | 4.88 | 5.01 | 4.77 | 4.09 | 4.49 | 5.31 |
| | 0.7 | 5.48 | 5.91 | 6.83 | 9.21 | 5.24 | 5.07 | 4.85 | 4.91 | 3.77 | 4.87 | 5.97 |
| $\mathbf{n}_3$ | 0.1 | 4.77 | 4.61 | 4.60 | 4.75 | 5.02 | 4.81 | 4.86 | 4.80 | 4.84 | 4.94 | 5.08 |
| | 0.3 | 5.14 | 5.11 | 4.73 | 4.99 | 5.09 | 5.19 | 4.76 | 5.11 | 4.56 | 4.74 | 4.95 |
| | 0.5 | 4.86 | 4.67 | 4.77 | 4.82 | 4.93 | 4.93 | 5.00 | 5.18 | 4.74 | 5.01 | 5.49 |
| | 0.7 | 4.69 | 4.78 | 4.43 | 4.88 | 5.18 | 4.80 | 4.64 | 4.98 | 4.57 | 5.17 | 5.37 |

Tables 3 and 4 display the empirical sizes and powers of $FHK_D, FHK_G, L^2_{nv}, L^2_{br}$ and $L^2_{rp}$ when $z_{ijr},\ r = 1, 2, \cdots, q;\ j = 1, 2, \cdots, n_i;\ i = 1, 2 \overset{i.i.d.}{\sim} N(0, 1)$. We may make the following conclusions:

- In terms of size-controlling, $FHK_D$ is conservative under both small and large sample sizes but it is liberal under medium sample sizes. This result shows that $FHK_D$ is unstable. In addition, choosing large $d$ may cause negative effects on the performance of $FHK_D$ — larger $d$ makes $FHK_D$ more conservative under small sample sizes and makes $FHK_D$ more liberal under medium sample sizes. When the parameter $\rho$ is large, the $L^2_{nv}$ has conservative empirical sizes while $L^2_{rp}$ has inflated empirical sizes. As the sample size increases, their empirical sizes become closer to the nominal size. Among the five tests, $FHK_G$ and $L^2_{br}$ slightly outperform others.

Table 4: Empirical powers (in percentages) of $FHK_D, FHK_G, L^2_{nv}, L^2_{br}$ and $L^2_{rp}$ for Simulation 2 when $z_{ijr},\ r=1,\cdots,q;\ j=1,\cdots,n_i;\ i=1,2 \overset{i.i.d.}{\sim} N(0,1)$.

| | | | $FHK_D$ | | | | $FHK_G$ | | | | $L^2_{nv}$ | $L^2_{br}$ | $L^2_{rp}$ |
|---|---|---|---|---|---|---|---|---|---|---|---|---|---|
| **n** | $\rho$ | $\omega$ | $d=1$ | $d=2$ | $d=3$ | $d=4$ | $d=1$ | $d=2$ | $d=3$ | $d=4$ | | | |
| $\mathbf{n}_1$ | 0.1 | 2.20 | 12.40 | 31.82 | 17.53 | 10.59 | 18.09 | 43.65 | 31.51 | 22.15 | 16.80 | 17.57 | 19.45 |
| | | 5.00 | 70.39 | 92.95 | 73.59 | 45.57 | 85.06 | 98.59 | 95.93 | 88.88 | 84.70 | 85.43 | 83.54 |
| | 0.3 | 1.00 | 9.26 | 20.54 | 11.78 | 7.86 | 13.54 | 27.42 | 19.89 | 15.11 | 14.20 | 15.40 | 17.92 |
| | | 2.80 | 71.63 | 91.55 | 70.71 | 43.23 | 85.68 | 97.89 | 94.43 | 87.68 | 86.46 | 87.31 | 86.31 |
| | 0.5 | 0.68 | 9.15 | 16.32 | 10.02 | 6.21 | 13.17 | 22.21 | 15.07 | 12.13 | 14.36 | 16.01 | 19.28 |
| | | 2.00 | 66.47 | 79.52 | 63.12 | 37.90 | 81.42 | 92.64 | 90.57 | 81.53 | 84.41 | 85.41 | 84.50 |
| | 0.7 | 0.55 | 9.85 | 13.01 | 8.82 | 5.86 | 14.40 | 18.01 | 13.96 | 11.02 | 13.21 | 16.45 | 19.43 |
| | | 1.60 | 64.74 | 56.38 | 50.84 | 32.82 | 79.55 | 78.21 | 80.54 | 73.80 | 80.88 | 83.06 | 83.30 |
| $\mathbf{n}_2$ | 0.1 | 2.00 | 23.04 | 48.27 | 33.20 | 28.90 | 17.57 | 53.69 | 39.49 | 30.60 | 15.06 | 15.72 | 16.13 |
| | | 4.30 | 87.67 | 98.52 | 92.36 | 81.11 | 85.05 | 99.47 | 97.97 | 94.99 | 85.01 | 85.45 | 82.79 |
| | 0.3 | 1.00 | 21.15 | 36.58 | 27.04 | 22.84 | 15.89 | 41.30 | 31.07 | 23.01 | 17.73 | 18.51 | 20.64 |
| | | 2.50 | 91.12 | 98.31 | 92.48 | 81.05 | 89.66 | 99.42 | 98.04 | 95.49 | 91.39 | 91.99 | 89.66 |
| | 0.5 | 0.50 | 13.25 | 17.31 | 13.84 | 14.21 | 9.83 | 18.67 | 14.93 | 11.14 | 11.32 | 12.31 | 13.82 |
| | | 1.50 | 77.30 | 86.74 | 74.70 | 60.85 | 72.82 | 90.17 | 86.34 | 77.82 | 77.28 | 78.39 | 76.48 |
| | 0.7 | 0.50 | 19.84 | 21.33 | 17.77 | 16.47 | 15.15 | 21.57 | 18.04 | 14.86 | 15.14 | 17.69 | 19.79 |
| | | 1.50 | 90.74 | 84.92 | 82.85 | 73.08 | 89.12 | 83.79 | 88.49 | 88.99 | 90.87 | 91.84 | 89.95 |
| $\mathbf{n}_3$ | 0.1 | 1.40 | 14.72 | 58.77 | 43.39 | 32.46 | 13.40 | 62.73 | 49.68 | 38.51 | 16.88 | 17.18 | 17.55 |
| | | 3.00 | 83.99 | 99.86 | 99.10 | 96.75 | 84.02 | 99.90 | 99.70 | 99.05 | 88.30 | 88.43 | 87.50 |
| | 0.3 | 0.70 | 14.21 | 44.55 | 32.54 | 24.21 | 12.97 | 48.36 | 37.56 | 29.93 | 23.52 | 24.02 | 25.43 |
| | | 1.40 | 70.47 | 97.77 | 92.59 | 83.64 | 69.49 | 98.52 | 95.81 | 91.30 | 84.54 | 84.77 | 84.60 |
| | 0.5 | 0.45 | 14.10 | 32.04 | 22.38 | 16.14 | 12.87 | 34.42 | 26.52 | 20.22 | 24.05 | 24.83 | 25.58 |
| | | 1.00 | 70.63 | 94.35 | 86.90 | 76.38 | 69.59 | 95.96 | 91.80 | 84.92 | 87.23 | 87.60 | 87.06 |
| | 0.7 | 0.35 | 16.21 | 26.11 | 19.15 | 14.78 | 15.02 | 27.92 | 22.11 | 17.61 | 22.95 | 24.12 | 24.63 |
| | | 0.80 | 74.39 | 84.87 | 82.11 | 70.25 | 73.62 | 86.61 | 87.17 | 80.47 | 85.65 | 86.21 | 85.87 |

- In most cases, $FHK_D$ and $FHK_G$ have largest empirical powers at $d = 2$ and their

Table 5: Empirical sizes (in percentages) of $FHK_D, FHK_G, L^2_{nv}, L^2_{br}$ and $L^2_{rp}$ for Simulation 2 when $z_{ijr},\ r=1,\cdots,q;\ j=1,\cdots,n_i;\ i=1,2 \overset{i.i.d.}{\sim} t_4/\sqrt{2}$.

| | | $FHK_D$ | | | | $FHK_G$ | | | | $L^2_{nv}$ | $L^2_{br}$ | $L^2_{rp}$ |
|---|---|---|---|---|---|---|---|---|---|---|---|---|
| **n** | $\rho$ | $d=1$ | $d=2$ | $d=3$ | $d=4$ | $d=1$ | $d=2$ | $d=3$ | $d=4$ | | | |
| $\mathbf{n}_1$ | 0.1 | 3.58 | 3.19 | 2.88 | 2.10 | 22.18 | 26.60 | 28.76 | 28.71 | 20.76 | 21.58 | 5.72 |
| | 0.3 | 3.56 | 3.03 | 2.90 | 2.27 | 23.18 | 27.21 | 27.95 | 29.31 | 21.82 | 23.09 | 6.29 |
| | 0.5 | 3.50 | 3.03 | 2.69 | 2.33 | 24.60 | 28.56 | 30.01 | 31.06 | 25.55 | 27.44 | 6.78 |
| | 0.7 | 3.06 | 2.97 | 2.31 | 1.94 | 28.29 | 32.61 | 32.77 | 33.07 | 30.22 | 34.68 | 7.19 |
| $\mathbf{n}_2$ | 0.1 | 4.46 | 4.16 | 5.11 | 6.21 | 23.84 | 29.85 | 33.83 | 34.52 | 22.13 | 22.71 | 5.69 |
| | 0.3 | 4.33 | 4.17 | 4.63 | 5.28 | 24.61 | 29.40 | 33.30 | 35.93 | 25.38 | 26.07 | 5.69 |
| | 0.5 | 4.33 | 4.44 | 4.80 | 5.67 | 26.02 | 32.49 | 34.63 | 36.28 | 30.39 | 31.77 | 6.32 |
| | 0.7 | 3.59 | 4.25 | 4.50 | 5.45 | 29.27 | 33.83 | 37.18 | 39.53 | 36.64 | 39.38 | 6.90 |
| $\mathbf{n}_3$ | 0.1 | 3.38 | 3.38 | 3.44 | 3.41 | 27.50 | 36.47 | 42.45 | 45.98 | 27.59 | 27.85 | 5.35 |
| | 0.3 | 3.85 | 3.38 | 3.30 | 3.51 | 29.24 | 37.11 | 42.28 | 46.47 | 31.37 | 31.62 | 5.19 |
| | 0.5 | 3.80 | 3.58 | 3.55 | 3.52 | 29.27 | 38.52 | 43.21 | 46.32 | 36.58 | 37.11 | 5.18 |
| | 0.7 | 3.55 | 3.27 | 2.98 | 3.04 | 32.63 | 39.72 | 44.74 | 48.60 | 46.08 | 47.33 | 5.36 |

empirical powers decrease with increasing the values of $d$ from 2 to 3, 4. This is not a surprise since $d=2$ is the correct number of functional principal components with the main differences of the two covariance functions located at the first two basis functions. This shows that the performances of $FHK_D$ and $FHK_G$ strongly depend on if a correct number of principal components is used. When the number of functional principal components is not well chosen, the powers of $FHK_D$ and $FHK_G$ may be much smaller than those of $L^2_{nv}, L^2_{br}$ and $L^2_{rp}$, as shown in the table.

- Unlike $FHK_D$ and $FHK_G$, $L^2_{nv}, L^2_{br}$ and $L^2_{rp}$ do not need to calculate the eigenvalues and choose the number of functional principal components. This could be a big advantage.

Tables 5 and 6 display the empirical sizes and powers of $FHK_D, FHK_G, L^2_{nv}, L^2_{br}$ and $L^2_{rp}$ when $z_{ijr},\ r=1,2,\cdots,q;\ j=1,2,\cdots,n_i;\ i=1,2 \overset{i.i.d.}{\sim} t_4/\sqrt{2}$. We may make the following conclusions:

- $FHK_G, L^2_{nv}$ and $L^2_{br}$ have too large empirical sizes and they do not work at all in this

Table 6: Empirical powers (in percentages) of $FHK_D, FHK_G, L^2_{nv}, L^2_{br}$ and $L^2_{rp}$ for Simulation 2 when $z_{ijr},\ r = 1, \cdots, q;\ j = 1, \cdots, n_i;\ i = 1, 2 \overset{i.i.d.}{\sim} t_4/\sqrt{2}$.

| | | | $FHK_D$ | | | | $FHK_G$ | | | | $L^2_{nv}$ | $L^2_{br}$ | $L^2_{rp}$ |
|---|---|---|---|---|---|---|---|---|---|---|---|---|---|
| **n** | $\rho$ | $\omega$ | $d=1$ | $d=2$ | $d=3$ | $d=4$ | $d=1$ | $d=2$ | $d=3$ | $d=4$ | | | |
| $\mathbf{n}_1$ | 0.1 | 2.60 | 12.60 | 36.43 | 18.81 | 10.09 | 37.14 | 71.50 | 64.30 | 58.56 | 36.26 | 37.17 | 23.84 |
| | | 6.50 | 61.89 | 88.71 | 64.83 | 36.93 | 90.19 | 99.70 | 99.35 | 98.55 | 89.97 | 90.33 | 82.79 |
| | 0.3 | 1.40 | 13.07 | 31.45 | 17.35 | 8.70 | 38.26 | 67.66 | 61.63 | 56.84 | 41.79 | 42.90 | 26.86 |
| | | 4.20 | 68.59 | 86.44 | 69.49 | 41.43 | 94.42 | 99.49 | 99.80 | 99.23 | 94.11 | 94.29 | 88.97 |
| | 0.5 | 1.00 | 11.56 | 23.99 | 14.84 | 7.76 | 37.99 | 61.04 | 57.89 | 53.63 | 45.65 | 47.87 | 27.78 |
| | | 2.80 | 61.02 | 66.97 | 60.52 | 35.22 | 90.28 | 96.41 | 98.71 | 97.89 | 91.73 | 92.24 | 84.82 |
| | 0.7 | 0.80 | 10.75 | 15.26 | 11.16 | 6.73 | 38.99 | 54.28 | 54.54 | 52.06 | 49.50 | 53.80 | 27.18 |
| | | 2.50 | 64.79 | 51.09 | 45.80 | 32.15 | 92.71 | 94.71 | 96.64 | 97.37 | 94.18 | 94.79 | 88.34 |
| $\mathbf{n}_2$ | 0.1 | 2.00 | 16.92 | 42.38 | 27.67 | 21.48 | 33.14 | 70.69 | 65.55 | 60.78 | 33.60 | 34.39 | 13.71 |
| | | 5.00 | 71.76 | 94.62 | 83.36 | 66.80 | 86.41 | 99.86 | 99.53 | 98.82 | 85.75 | 86.10 | 68.61 |
| | 0.3 | 1.00 | 16.36 | 33.07 | 21.52 | 16.87 | 33.77 | 61.57 | 58.33 | 57.13 | 38.74 | 39.81 | 14.94 |
| | | 4.00 | 87.45 | 95.44 | 91.60 | 79.29 | 97.03 | 99.75 | 99.99 | 99.98 | 97.41 | 97.49 | 88.25 |
| | 0.5 | 1.00 | 27.97 | 43.25 | 33.58 | 25.24 | 46.37 | 72.81 | 72.33 | 69.27 | 58.17 | 59.69 | 27.54 |
| | | 2.80 | 85.20 | 86.54 | 87.32 | 76.07 | 95.79 | 98.43 | 99.79 | 99.84 | 96.61 | 96.71 | 86.59 |
| | 0.7 | 0.80 | 27.18 | 31.70 | 27.99 | 21.89 | 47.02 | 64.01 | 68.28 | 66.06 | 62.65 | 65.33 | 25.49 |
| | | 2.10 | 80.60 | 73.02 | 72.22 | 66.24 | 93.13 | 94.62 | 96.71 | 98.37 | 95.78 | 96.04 | 82.56 |
| $\mathbf{n}_3$ | 0.1 | 1.20 | 7.44 | 43.31 | 29.03 | 19.60 | 31.67 | 71.80 | 68.35 | 67.97 | 34.90 | 35.11 | 9.38 |
| | | 2.60 | 40.99 | 96.16 | 88.81 | 78.12 | 64.93 | 99.49 | 98.95 | 97.66 | 69.46 | 69.69 | 45.66 |
| | 0.3 | 0.70 | 10.00 | 42.72 | 29.89 | 20.67 | 33.95 | 72.06 | 69.87 | 69.13 | 48.15 | 48.48 | 14.29 |
| | | 1.50 | 47.90 | 95.92 | 88.95 | 77.02 | 71.29 | 99.29 | 98.55 | 97.90 | 84.94 | 85.20 | 58.34 |
| | 0.5 | 0.50 | 11.41 | 35.94 | 25.40 | 18.05 | 35.96 | 67.55 | 66.38 | 66.98 | 58.66 | 59.45 | 17.03 |
| | | 1.00 | 42.60 | 86.03 | 78.97 | 65.20 | 67.38 | 96.64 | 96.12 | 94.90 | 88.59 | 88.98 | 56.56 |
| | 0.7 | 0.45 | 15.47 | 32.57 | 27.46 | 19.90 | 41.53 | 67.20 | 69.77 | 69.86 | 69.82 | 71.11 | 19.36 |
| | | 0.85 | 46.22 | 69.18 | 76.47 | 65.17 | 71.41 | 89.92 | 95.51 | 95.25 | 93.49 | 93.81 | 61.60 |

simulation setting. This is expected since they are developed only for Gaussian functional data.

- In terms of sizes, $FHK_D$ is rather conservative even when the correct number of functional principal components, $d = 2$, is used, especially under the small and large sample sizes, while $L^2_{rp}$ works reasonably well under the large sample size, and is slightly inflated under the small and medium sample sizes when $\rho$ is large. In terms of powers, $FHK_D$ generally outperforms $L^2_{rp}$, except some cases when sample sizes are not large and when the data are less correlated.

Table 7: Empirical sizes and powers (in percentages) of $FHK_D, FHK_G, L^2_{nv}, L^2_{br}$ and $L^2_{rp}$ for Simulation 2 with $\mathbf{n}_3 = [75, 85]$ when $z_{ijr},\ r = 1, \cdots, q;\ j = 1, \cdots, n_i;\ i = 1, 2$ are i.i.d. $N(0, 1)$ under the new simulation scheme.

| | | $FHK_D$ | | | | $FHK_G$ | | | | $L^2_{nv}$ | $L^2_{br}$ | $L^2_{rp}$ |
|---|---|---|---|---|---|---|---|---|---|---|---|---|
| $\rho$ | $\omega$ | $d=1$ | $d=2$ | $d=3$ | $d=4$ | $d=1$ | $d=2$ | $d=3$ | $d=4$ | | | |
| | 0.00 | 5.34 | 4.48 | 4.35 | 4.46 | 5.34 | 4.76 | 4.76 | 4.83 | 4.36 | 4.45 | 4.92 |
| 0.1 | 0.64 | 4.77 | 4.56 | 4.72 | 4.50 | 4.97 | 4.88 | 4.77 | 4.27 | 53.20 | 54.14 | 56.83 |
| | 0.76 | 4.87 | 4.86 | 5.07 | 4.57 | 4.96 | 5.01 | 5.38 | 4.54 | 94.46 | 94.83 | 94.36 |
| | 0.00 | 4.99 | 4.52 | 4.64 | 4.23 | 5.00 | 4.66 | 5.23 | 4.40 | 4.42 | 4.58 | 4.82 |
| 0.3 | 0.64 | 4.80 | 4.79 | 4.40 | 4.71 | 5.01 | 4.86 | 4.70 | 4.84 | 50.61 | 52.08 | 54.35 |
| | 0.77 | 5.11 | 5.05 | 4.51 | 5.15 | 5.20 | 5.05 | 4.74 | 5.05 | 94.86 | 95.18 | 94.77 |
| | 0.00 | 5.22 | 4.69 | 4.50 | 4.80 | 5.07 | 4.91 | 4.63 | 5.05 | 4.90 | 5.07 | 5.53 |
| 0.5 | 0.66 | 4.88 | 4.81 | 4.71 | 4.70 | 5.05 | 4.90 | 4.98 | 5.20 | 53.23 | 54.94 | 56.49 |
| | 0.78 | 4.79 | 4.99 | 4.32 | 4.48 | 4.92 | 5.16 | 4.80 | 4.69 | 93.27 | 94.00 | 93.34 |
| | 0.00 | 5.01 | 5.45 | 4.97 | 4.48 | 5.01 | 5.49 | 5.37 | 4.94 | 4.55 | 5.00 | 5.32 |
| 0.7 | 0.69 | 4.81 | 5.11 | 4.35 | 4.66 | 4.90 | 5.30 | 4.62 | 4.87 | 51.20 | 54.44 | 54.97 |
| | 0.82 | 4.96 | 4.78 | 4.70 | 4.85 | 4.99 | 4.70 | 5.21 | 5.11 | 92.02 | 93.10 | 92.78 |

In some situations, however, $L^2_{nv}, L^2_{br}$ and $L^2_{rp}$ can have much higher powers than $FHK_D$ and $FHK_G$. We can show this via making a minor change of the previous simulation scheme. We continue to use the data generating model (4.1) but we set $\eta_i(t) = 0, i = 1, 2, \cdots, k$ for simplicity and only present the large sample size result for space saving. In addition, we increase the number of basis functions to $q = 25$ and set $\lambda_{1r} = \rho^{r-1}, r = 1, 2, \cdots, q$, $\lambda_{2r} = \rho^{r-1}, r = 1, 2, \cdots, q-1$ and $\sqrt{\lambda_{2q}} = \sqrt{\lambda_{1q}} + \omega$ so that the differences of the covariance functions of the functional samples are located in the space spanned by the last eigenfunction. Since the information is located in high frequencies, $FHK_D$ and $FHK_G$ will be less powerful in detecting the differences of the covariance functions. This is not the case for $L^2_{nv}, L^2_{br}$ and $L^2_{rp}$ since these $L^2$-norm based tests use all the information provided by the data. The simulation results presented in Table 7 indeed show this is true. From this table, we can see that the powers of $FHK_D$ and $FHK_G$ are about the same with increasing the values of $\omega$ but powers of $L^2_{nv}, L^2_{br}$ and $L^2_{rp}$ become larger as $\omega$ increases.

# 5 Applications to Two Real Data Examples

In this section, we shall present the applications of $L^2_{nv}, L^2_{br}, L^2_{rp}$ to two real data examples. Throughout this section, the p-values of $L^2_{rp}$ were obtained based on 10000 runs of random permutations.

## 5.1 The Medfly Data

In this subsection, we present some applications of $L^2_{nv}, L^2_{br}, L^2_{rp}$ to check if the cell covariance functions of the medfly data are the same. The medfly data, which recorded daily egg-laying numbers of 1000 medflies (Mediterranean fruit flies), have been analyzed by several authors in the literature, including Müller and Stadtmüller (2005) and Fremdt et al. (2013) among others. Thanks to Professor Hans-Georg Müller and Professor Carey's laboratory, this medfly data set is available at http://anson.ucdavis.edu/~mueller/data/data.html. Previous studies indicate that the fecundity may be associated with the individual mortality and longevity.

We picked up 534 medflies who lived at least 34 days and studied both the absolute and relative counts of eggs laid by the 534 medflies in the first 30 days. A relative count is defined as the ratio of absolute count in each day to the total number of eggs laid by each medfly. These medflies are classified into two groups: long-lived and short-lived. The long-lived group includes 278 medflies who lived 44 days or longer and the short-lived group includes 256 medflies who lived less than 44 days.

Fremdt et al. (2013) has considered testing if the covariance functions of the long-lived group and the short-lived group are the same. We can also apply $L^2_{nv}, L^2_{br}, L^2_{rp}$ for this problem. Actually, based on the absolute counts, the test statistic $T_n = 2.9774e8$ and the p-values of $L^2_{nv}, L^2_{br}, L^2_{rp}$ are $0.3017, 0.2999$ and $0.1228$ respectively. These p-values show that there is no strong evidence against the null hypothesis that the covariance functions of the long-lived group and the short-lived group are the same. This conclusion is consistent with the one made by the $FHK_D$ test described in Fremdt et al. (2013). Based on the relative counts, on the other hand, the associated test statistic $T_n = 0.0191$ and the p-values of $L^2_{nv}, L^2_{br}, L^2_{rp}$ are now $0, 0$ and $0.0025$ respectively. These p-values show that there is very strong evidence against the null hypothesis. This conclusion is again consistent with the one made by the $FHK_D$ test. However,the $FHK_G$ test can not obtain the right result. Besides, choosing different $d$ may get different p-values.

According to Fremdt et al. (2013), both the absolute counts and the relative counts have a strong deviation from normality which can be easily verified by QQ-plots. Therefore, the p-values of $L^2_{rp}$ will be more reliable than those of $L^2_{nv}$ and $L^2_{br}$. Although $L^2_{nv}$ and $L^2_{br}$ are based on the Gaussian assumption, these two tests give consistent result for both the absolute counts and relative counts while $FHK_G$ may give a misleading conclusion because the result of $FHK_G$ varies depending on the selection of empirical functional principal components as shown in Fremdt et al. (2013).

It is also possible to classify the medflies into three groups. The first group consists of the long-lived medflies who lived 50 days or longer, the second group consists of the medium-lived medflies who lived at least 40 days but no longer than 50 days, and the third group consists of the short-lived medflies who lived less than 40 days. Of the 534 medflies, 180 are long-lived, 180 are medium-lived and 174 are short-lived. Of interest is to test if the covariance functions of the three groups of medflies are the same.

Based on the absolute counts, the associated test statistic $T_n = 5.7069e8$ and the p-values of $L^2_{nv}, L^2_{br}, L^2_{rp}$ are $0.3132, 0.3107$ and $0.1030$ respectively. Thus, again, there is no strong evidence against that the covariance functions of the three groups are the same. Based on the relative counts, the associated test statistic $T_n = 0.0337$ and the p-values of $L^2_{nv}, L^2_{br}, L^2_{rp}$ are $0, 0$, and $0.0123$ respectively. Thus, again, there is strong evidence against that the covariance functions of the three groups are the same. These conclusions are consistent with those obtained based on the comparison of the covariance functions of the long-lived and short-lived medflies described above.

### 5.2 The orthosis data

In this subsection, we present some applications of $L^2_{nv}, L^2_{br}, L^2_{rp}$ to check if the cell covariance functions of the orthosis data are the same. The orthosis data set was kindly provided by Dr. Brani Vidakovic via email communication. It has been previously studied by a number of authors, including Abramovich et al. (2004), Abramovich and Angelini (2006), Antoniadis and Sapatinas (2007), and Cuesta-Albertos and Febrero-Bande (2010) among others.

To better understand how muscle copes with an external perturbation, the orthosis data were acquired and computed by Dr. Amarantini David and Dr. Martin Luc (Laboratoire Sport et Performance Motrice, EA 597, UFRAPS, Grenoble University, France). The data set recorded the moments at the knee of 7 volunteers under 4 experimental conditions (control, orthosis, spring 1, spring 2), each 10 times at equally spaced 256 time points. Figure 1 displays the raw curves of the orthosis data set, with each panel showing 10 raw curves. Figure 2 shows the 4 estimated cell covariance functions for the fifth volunteer under all the 4 conditions. Based on these two figures, it seems that the cell covariance functions are not exactly the same.

Figure 1: Raw curves of the orthosis data.

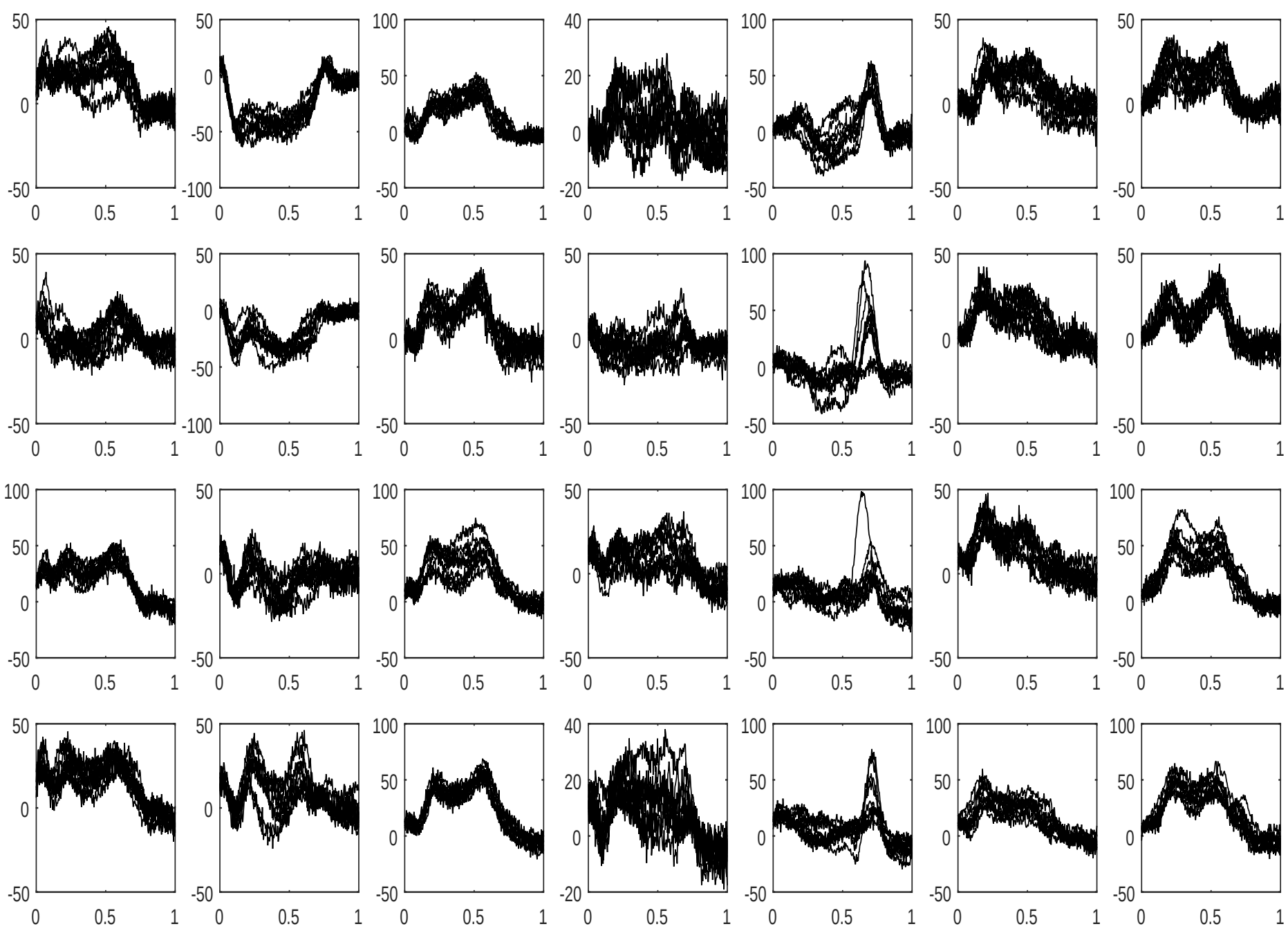

Figure 2: Estimated cell covariance functions of the orthosis data for the fifth volunteer under 4 treatment conditions.

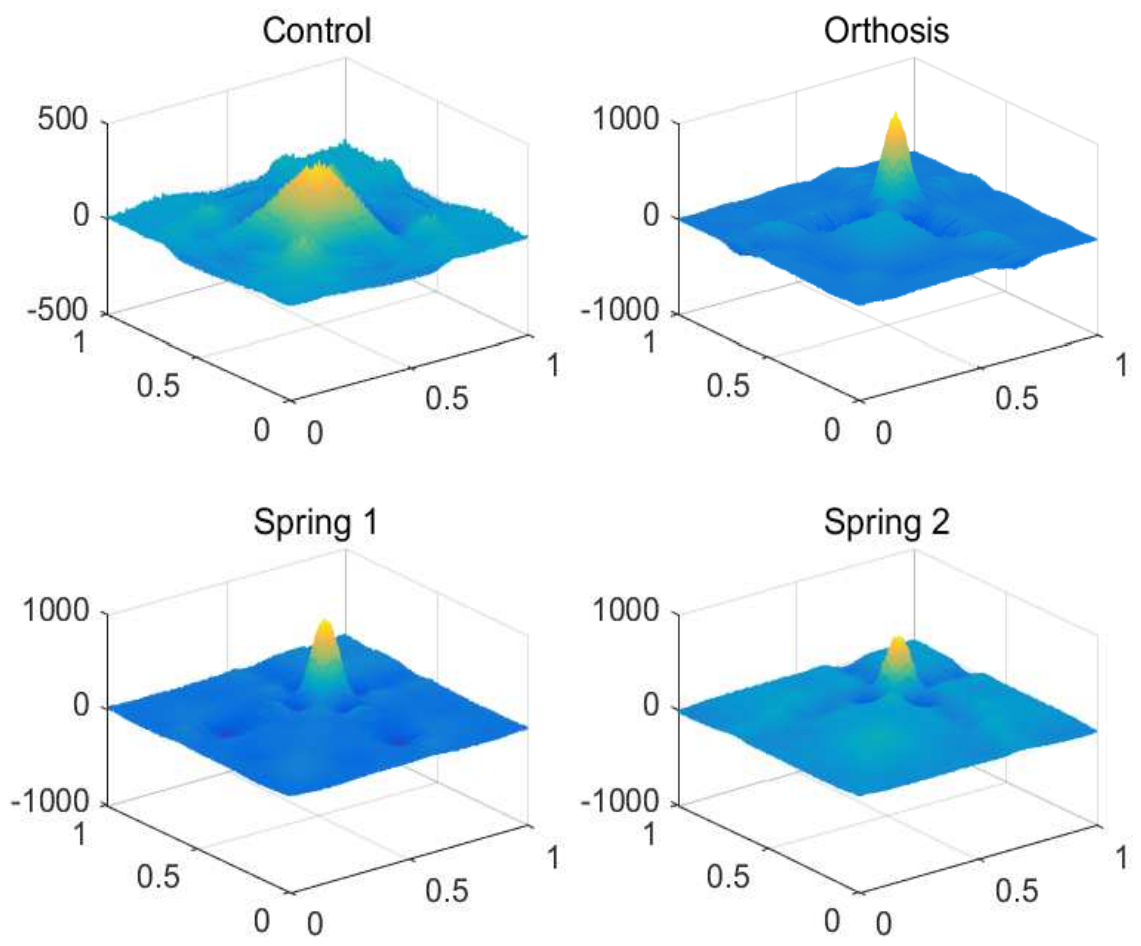

We first applied $L^2_{nv}, L^2_{br}, L^2_{rp}$ to test if all the 28 cell covariance functions are the same. It is easy to obtain that the test statistic $T_n = 1.5661e10$ using (2.6). To apply $L^2_{nv}$, by (2.10) and (2.11), we obtained $\mathrm{tr}(\hat{\varpi}) = 2.9198e8$, $\mathrm{tr}(\hat{\varpi}^{\otimes 2}) = 5.1118e15$ so that $\hat{\beta} = 1.7507e7, \hat{d} = 450.29$. The resulting p-value of $L^2_{nv}$ is then 0. To apply $L^2_{br}$, we obtained $\widehat{\mathrm{tr}(\varpi)} = 2.9051e8$, $\widehat{\mathrm{tr}(\varpi^{\otimes 2})} = 4.9242e15$ using (2.13) and then $\hat{\beta} = 1.6950e7, \hat{d} = 462.75$ using (2.12). The resulting p-value of $L^2_{br}$ is also 0. Similarly, the resulting p-value of $L^2_{rp}$ is again 0. These p-values demonstrate that the 28 cell covariance functions of the orthosis data are unlikely to be the same.

Second, we applied $L^2_{nv}, L^2_{br}, L^2_{rp}$ to check if the 4 estimated cell covariance functions for the fifth volunteer under the 4 different conditions are the same. The resulting test statistic is $T_n = 5.5510e9$ and the resulting p-values of $L^2_{nv}, L^2_{br}, L^2_{rp}$ are $0.0040, 0.0011$ and $0.0150$ respectively. These resulting p-values show that the 4 cell covariance functions under consideration are unlikely to be the same which is consistent with what we observed from Figure 2.

Finally, we applied $L^2_{nv}, L^2_{br}, L^2_{rp}$ to check if the 4 estimated cell covariance functions for the first volunteer under the 4 different conditions are the same. The resulting test statistic is $T_n = 5.7008e8$ and the resulting p-values of $L^2_{nv}, L^2_{br}, L^2_{rp}$ are $0.4670, 0.4050$ and $0.2076$ respectively, showing that there is no strong evidence that the 4 estimated cell covariance functions for the first volunteer under the 4 different conditions are the same.

## Appendix

Technical proofs and additional contents are available in supplementary materials.